\documentclass[letterpaper, 12 pt]{article}
\usepackage{amsthm, amsmath, amssymb, amsfonts, graphics}
\usepackage[all]{xy}

\newcommand{\setleftmargin}[1]{
        \addtolength{\textwidth}{\oddsidemargin}
        \addtolength{\textwidth}{1in}
        \addtolength{\textwidth}{-#1}
        \setlength{\oddsidemargin}{-1in}
        \addtolength{\oddsidemargin}{#1}
        \setlength{\evensidemargin}{\oddsidemargin}
}
\newcommand{\setrightmargin}[1]{
        \setlength{\textwidth}{8.5in}
        \addtolength{\textwidth}{-\oddsidemargin}
        \addtolength{\textwidth}{-1in}
        \addtolength{\textwidth}{-#1}
}
\setleftmargin{1 in}
\setrightmargin{1 in}

\def\ZZ{\mathbb{Z}}
\def\QQ{\mathbb{Q}}
\def\RR{\mathbb{R}}
\def\CC{\mathbb{C}}
\def\PP{\mathbb{P}}

\def\a{\boldsymbol{\alpha}}
\def\b{\boldsymbol{\beta}}
\def\c{\boldsymbol{\gamma}}

\def\Trop{\mathcal{T}}

\def\TRR{\RR \{ t \}}
\def\TCC{\CC \{ t \}}
\def\TRRp{{\TRR_{+}}}

\def\Id{\mathrm{Id}}

\def\Mat{\mathrm{Mat}}
\def\Her{\mathrm{Her}}
\def\GL{\mathrm{GL}}
\def\PD{\mathrm{PD}}

\def\Eigen{\mathrm{Eigen}}
\def\Crit{\mathrm{Sing}}
\def\Sing{\mathrm{Sing}}
\def\Root{\mathrm{Root}}

\def\HIVE{\mathtt{HIVE}}
\def\HH{\mathcal{H}}
\def\OO{\mathcal{O}}
\def\One{\boldsymbol{1}}
\def\C{\mathcal{C}}
\def\D{\mathcal{D}}
\def\P{\mathcal{P}}
\def\Tau{\textbf{T}}
\def\A{\mathcal{A}}
\def\sym{\sim}

\def\dec{\mathrm{dec}}
\def\add{\mathrm{add}}
\def\mult{\mathrm{mult}}

\def\isomorph{\cong}

\newtheorem{conj}{Conjecture}
\newtheorem{Theorem}[conj]{Theorem}
\newtheorem{prop}[conj]{Proposition}
\newtheorem{Lemma}[conj]{Lemma}

\newtheorem{cor}[conj]{Corollary}
\newtheorem{question}[conj]{Question}

\newtheorem*{MT}{Main Theorem}
\newtheorem*{HT}{The Hive Theorem}
\newtheorem*{KT}{Klyachko's Theorem}
\newtheorem*{Tar}{Tarski's Theorem}
\newtheorem*{VT}{Viro's Patchworking Theorem}
\newtheorem*{VinCrit}{Vinnikov's Criterion}
\newtheorem*{NC}{A New Criterion for Horn's Problem}

\begin{document}

\title{Horn's Problem, Vinnikov Curves and the Hive Cone}
\author{David E Speyer \\ \texttt{speyer@math.berkeley.edu}}
\maketitle

\abstract{A Vinnikov curve is a projective plane curve which can be written in the form $\det(xX+yY+zZ)=0$ for $X$, $Y$ and $Z$ positive definite Hermitian $n \times n$ matrices. Given three $n$-tuples of positive real numbers, $\a$, $\b$ and $\c$, there exist $A$, $B$ and $C \in \GL_n \CC$ with singular values $\a$, $\b$ and $\c$ and $ABC=1$ if and only if there is a Vinnikov curve passing through the $3n$ points $(-1: \alpha_i^2:0)$, $(0:-1:\beta_i^2)$ and $(\gamma_i^2:0:-1)$. Knutson and Tao proved that another equivalent condition for such $A$, $B$ and $C$ to exist is that there is a hive (defined within) whose boundary is $(\log \alpha, \log \beta, \log \gamma)$. The logarithms of the coefficients of $F$ approximately form such a hive; this leads to a new proof of Knutson and Tao's result. This paper uses no representation theory and essentially no symplectic geometry. In their place, it uses Viro's patchworking method and a topological description of Vinnikov curves.}

\section{Introduction}

This paper considers two problems of real algebraic geometry. The first, known as Horn's problem, asks about the eigenvalues of three Hermitian matrices with sum zero. Horn's problem has connections to representation theory of the general linear group and modules over a discrete valuation ring, as well as numerous other extremization problems in the theory of matrices. See \cite{Ful} for a survey of these problems and their relationships and \cite{Hor} for Horn's original paper. Knutson and Tao have given an elegant solution to Horn's problem (described later in this paper) in terms of combinatorial objects known as hives. The apperarance of hives is somewhat unexplained and the proofs of the Knutson-Tao characterization either go through representation theory or symplectic geometry.

We will show that Horn's problem is related to another problem, first raised by Lax and investigated extensively by Vinnikov, of determining which polynomials $F(x,y,z)$ can be written as $\det(xX+yY+zZ)$ with $X$, $Y$ and $Z$ positive definite. This problem has relevance in control theory. In \cite{Vin2}, Vinnikov gave a characterization of these polynomials in terms of the topology of their real points. 

We now state our result, leaving some of the details of our notation to be defined in the succeeding sections. Let $\TCC$ denote the field of Laurent series in fractional powers of $t^{-1}$, also known as Pusieux series. (We will impose the condition that our power series converge for $t$ sufficiently large in order to simplify a technical issue in the proof. The results stated here are true whether or not we impose convergence.) We write $l: \TCC^* \to \QQ$ for the map which takes a Pusieux series to its leading exponent. Let $\TRR$ denote the subfield of Pusieux series with real coefficients and $\TRRp$ the subsemifield of $\TRR$ consisting of those power series with positive leading term and real coefficients.

A Vinnikov polynomial over $\TRR$ is a homogeneous degree $n$ polynomial $F(x,y,z) \in \TRR[x,y,z]$ which can be written as $\det(xX(t)+yY(t)+zZ(t))$ with $X$, $Y$ and $Z$ positive definite matrices with entries in $\TCC$. We will show in Lemma~\ref{PosCoeff} that the coefficients of a Vinnikov polynomial lie in $\TRRp$.

\begin{MT}
Let $F = \sum F_{ijk}(t) x^i y^j z^k$ be a homogenous degree $n$ polynomial with coefficients in $\TRRp$ and set $h_{ijk}=l(F_{ijk}) $. If $F$ is a Vinnikov polynomial, $h$ is a hive and, if $h$ is a strict hive, $F$ is a Vinnikov polynomial. (See section~\ref{Hives} for definitions of ``hive'' and ``strict hive''.)
\end{MT}

In the rest of this section, we will describe needed vocabulary and concepts. Our basic notation concerning matrices is as follows: 

Let $\Mat_n$ denote the space of $n$ dimensional complex matrices. $\GL_n$ is the group of invertible matrices over $\CC$ of order $n$. $U_n$ is the unitary group over $\CC$ of order $n$. $\Her_n$ is the vector space of Hermitian (not necessarily non-degenerate) matrices over $\CC$ of order $n$. $\PD_n$ is the set of positive definite Hermitian matrices over $\CC$ of order $n$.

We will denote the conjugate transpose of $A$ by $A^*$. We will almost always drop the subscript $n$. We will sometimes write $\Mat_n(\CC)$ and so forth to emphasize that our matrices have entries in $\CC$ and not some other field. If $F$ is a polynomial, we denote its (complex) zero locus by $Z(F)$. 

\subsection{Horn's Problem}

Fix an integer $n$. Let $\a$, $\b$ and $\c \in \RR^n$. Horn's additive problem asks whether there exists a triple $(A,B,C)$ of $n \times n$ Hermitian matrices with $A+B+C=0$ and eigenvalues $\a$, $\b$ and $\c$. If so, we will say that \emph{Horn's additive problem is solvable for $\a$, $\b$ and $\c$}. 

Let $\OO^{\add}$ denote the space of triples of Hermitian matrices with sum $0$, modulo the action of the unitary group by conjugation. Write $(\RR)^n_{\dec}$ for the space of weakly decreasing $n$-tuples of real numbers and let $\Eigen : \OO^{\add} \to ((\RR^n)_\dec)^3$ denote the map which takes $(A,B,C)$ to their eigenvalues. Horn's additive problem can therefore be described as asking to characterize the image of $\Eigen$.

The \emph{singular values} of a matrix $A \in \Mat_n \CC$ are defined to be the positive square roots of the eigenvalues of $A^* A$. Horn's multiplicative problem asks, if $\a$, $\b$ and $\c$ are in $(\RR_+)^n$, whether there exists a triple $(A,B,C)$ of invertible matrices with $ABC=1$ and singular values $\a$, $\b$ and $\c$. If so, we say \emph{Horn's multliplicative problem is solvable for $(\a, \b, \c)$}.

Set $\OO^{\mult}$ to denote the set of triples $(A,B,C)$ of invertible matrices with $ABC=1$ modulo the action of  $U^3$, where the action of $U^3$ is 
$$(U,V,W) : (A,B,C) \mapsto (W^{-1} A U, U^{-1} B V, V^{-1} C W).$$
Denote by $\Crit$ the map which sends $\OO^{\mult}$ to $((\RR_+^n)_{\dec})^3$ by sending $(A,B,C)$ to their singular values.

\begin{KT}
Horn's additive problem is solvable for $(\a, \b, \c)$ if and only if Horn's multiplicative problem is solvable for $(e^{\a}, e^{\b}, e^{\c})$. In fact, there is a diffeomorphism $\alpha: \OO^{\add} \to \OO^{\mult}$ such that $\Crit \circ \alpha=\exp \circ \Eigen$. 
\end{KT}

\begin{proof}
The first sentence is due to \cite{Kly}. The second is essentially due to \cite{AMW}. What \cite{AMW} actually do is to construct diffeomorphisms between each fiber of $\Eigen$ and the corresponding fiber of $\Sing$, but one can check that they do glue together. To see this carried out and put in a more general context, see \cite{HE}.
\end{proof}

\subsection{Hives and Related Notions}\label{Hives}

Let $\Delta_n$ (usually abreviated $\Delta$) denote the set of triples $(i,j,k)$ of nonnegative integers such that $i+j+k=n$. We will visualize the points of $\Delta$ as arranged in a triangle. For any set $S$, the symbol $S^{\Delta}$ denotes the collection of functions $\Delta \to S$. 

A function $h_{ijk} \in \RR^{\Delta}$ is called a \emph{hive} if it obeys the inequalities
\begin{eqnarray*}
h_{(i+1)j(k-1)}+h_{i(j+1)(k-1)} & \geq & h_{ijk}+h_{(i+1)(j+1)(k-2)} \\
h_{(i+1)(j-1)k}+h_{i(j-1)(k+1)} & \geq & h_{ijk}+h_{(i+1)(j-2)(k+1)} \\
h_{(i-1)(j+1)k}+h_{(i-1)j(k+1)} & \geq & h_{ijk}+h_{(i-2)(j+1)(k+1)}
\end{eqnarray*}

It is called a \emph{strict hive} if all of these inequalities hold strictly. The set of hives forms a polyhedral cone, which we will denote by $\HIVE_n \subset \RR^{\Delta}$ and will use the same notation when considering it as a cone in the quotient $\RR^{\Delta}/\RR\One$.

The boundary of a hive $h_{ijk}$, denoted $\partial_0 h_{ijk}$, is the vector in $((\RR^n)_{\dec})^3$ given by 
$$
\begin{array}{c}
((h_{(n-1)10} - h_{n00}, h_{(n-2)20}-h_{(n-1)10}, \ldots, h_{0n0}-h_{1(n-1)0}), \phantom{)} \\
\phantom{(} (h_{0(n-1)1} - h_{0n0}, h_{0(n-2)2}-h_{0(n-1)1}, \ldots, h_{00n}-h_{01(n-1)}), \phantom{)} \\ 
\phantom{(} (h_{10(n-1)} - h_{00n}, h_{20(n-2)}-h_{10(n-1)}, \ldots, h_{n00}-h_{(n-1)01})).
\end{array}
$$

For a more geometric description of hives, see proposition~\ref{HivesandTriangulations}.

We can now state the Knutson-Tao solution to Horn's problem.

\begin{HT} 
Horn's additive problem is solvable for $(\a, \b, \c)$ if and only if $(\a, \b, \c)$ are the boundary of a hive.
\end{HT}

\begin{proof}
See the appendix of \cite{KTW}.
\end{proof}

One of our goals in this paper will be to give a new proof of this result.

\textbf{Remark:}
The result conjectured in Horn's original paper, \cite{Hor}, and which is proved in the Knutson and Tao papers, is a recursive description of a set of inequalities on $(\a, \b, \c)$ that are necessary and sufficent for Horn's additive problem to be solvable for $(\a, \b, \c)$. It is not obvious that these inequalties are equivalent to being the boundary of a hive, but it is true. The claim that these inequalities do in fact characterize those $(\a, \b, \c)$ for which Horn's problem is solvable is called Horn's conjecture.

\subsection{Vinnikov Curves}

We will define a \emph{Vinnikov curve} to be a projective plane curve given by a homogeneous degree $n$ polynomial $F(x,y,z)$ of the form $\det(xX+yY+zZ)$ with $X$, $Y$ and $Z$ positive definite hermitian matrices. These are also known as hyperbolic curves. We will denote the set of such polynomials by $\HH_n \in \RR[x,y,z]_n$. (The subscript $n$ denotes that we are dealing with homogenous polynomials of degree $n$.) Often we will work in $\RR[x,y,z]_n/\RR^{\times}$, thus equating proportional polynomials that cut out the same curve, and will also denote the image of $\HH_n$ in this quotient by $\HH_n$.

Set $\P_n$ to denote the set of triples $(X,Y,Z)$ of positive definite matrices modulo the action of  $\GL$, where the action of $\GL_n$ is 
$$G : (X,Y,Z) \mapsto (G^* X G, G^* Y G, G^* Z G).$$
Denote by $\delta$ the map which sends $\P$ to $\RR[x,y,z]_n/\RR^{\times}$ by sending $(X,Y,Z)$ to $\det(xX+yY+zZ)$, so $\HH$ is the image of $\delta$.

Vinnikov has made a detailed study of the map $\delta$. In this section, we summarize the results we need.

\begin{Lemma}\label{PosCoeff}
The coefficients of a polynomial in $\HH_n$ are positive. More generally, let $X_1$, \dots, $X_r \in \PD$ and set 
$$F(x_1, \ldots, x_r)=\det \left( \sum_{i=1}^r x_i X_i \right).$$
Then the coefficients of $F$ are positive real numbers.
\end{Lemma}

\begin{proof}
Our proof is by induction on $r$. If $r=1$, this says that the determinant of a positive definite matrix is positive.

We can choose an orthonormal basis with respect to $X_r$, so without loss of generality $X_r=\Id$. For $M \in \Mat$ and $I \subseteq [n]$, let $M_I$ denote the sub-matrix with rows and columns drawn from $I$.

We have
$$F=\det (x_r\Id+ \sum_{i=1}^{r-1} x_i X_i )=\sum_{I \subseteq [n]} x_r^{|[n] \setminus I|} \det ( \sum_{i=1}^{r-1} x_i (X_i)_I ).$$
$(X_i)_I$ is the restriction of the positive definite form $X_i$ to a subspace, so it is positive definite. Thus, by induction, every term in this sum has positive coefficients, so $F$ does as well.
\end{proof}

\begin{Theorem}\label{AnyCurve}
Every smooth $F \in \CC[x,y,z]_n$ can be written as $\delta(X,Y,Z)$ for $(X,Y,Z) \in \Mat_n(\CC)$.
\end{Theorem}

\begin{proof}
See \cite{Vin1}.
\end{proof}

\textbf{Remark:} The analogous statement for poynomials in more than three variables is false.

We now prove some lemmas that will allow us to reduce various problems to the case of smooth curves:

\begin{Lemma}\label{SmoothareDense}
The smooth curves are dense in $\HH$.
\end{Lemma}

\begin{proof}
Consider the map $\delta: \Mat_n(\CC)^3 \to \CC[x,y,z]_n$ via $(X,Y,Z) \mapsto \det(Xx+Yy+Zz)$. By Theorem~\ref{AnyCurve}, every smooth curve is in the image of $\delta$, so $ \delta$ is a dominant map. Let $\D \subset \CC[x,y,z]_n$ be the subset of non-smooth curves; we have that $\delta^{-1}(\D)$ is a proper subvariety of $\Mat(\CC)^3$.

Now, let $(X, Y, Z) \in \PD^3$. We must show that we can perturb $(X,Y,Z)$ to $(X',Y',Z') \in \PD^3$ such that $\delta(X',Y',Z') \not\in \D$. 

\begin{Lemma}
Let $V \subset \CC^N$ be a proper algebraic subvariety and let $v \in \CC^N$. Then any arbitrary neighborhood of $v$ contains points in $(v+\RR^N) \setminus V$.
\end{Lemma}

\begin{proof}
Suppose this were false. Without loss of generality, $v=0$. Let $f$ be any polynomial vanishing on $V$. Then $f$ is zero on a neighborhood of the identity in $\RR^N$ and hence $f=0$. But, if the only polynomial vanishing on $V$ is the zero polynomial, then $V=\CC^N$.
\end{proof}

Applying this lemma with $\CC^N=\Mat_n^3=\Her_n^3 \otimes_{\RR} \CC$, $\RR^N=\Her_n^3$, $V=\D$ and $v=(X,Y,Z)$, we see that there are points $(X',Y',Z')$ of $\Her^3$ arbitrarily close to $(X,Y,Z)$ with $\delta(X',Y',Z') \not \in D$. As $\PD^3$ is open in $\Her^3$, there will be points in $\PD^3 \setminus \delta^{-1}(\D)$ arbitrarily close to $(X,Y,Z)$.
\end{proof}

\begin{Lemma} \label{proper}
The map $\delta$ is proper as a map to the subset of $\RR[x,y,z]_n/\RR^{\times}$ where $F(1,0,0)$, $F(0,1,0)$ and $F(0,0,1)$ are nonzero. (Proper, in this context, means that the preimage of a compact set is compact.) In particular, $\HH_n$ is closed in this space.
\end{Lemma}

\begin{proof}
Consider some compact set in $\RR[x,y,z]_n/\RR^{\times}$ where $F(1,0,0)$ is always nonzero, we may choose our representative modulo the $\RR^{\times}$ action such that $F(1,0,0)$ is always $1$. We may furthermore take $X$ to always be $\Id$. Then the eigenvalues of $Y$ and $Z$ are the roots of $F(-u,1,0)$ and $F(-u,0,1)$. As the coefficients of $F$ are bounded, these eigenvalues remain in a bounded portion of $(\RR^{\times})^n$. A set of Hermitian matrices with bounded eigenvalues is itself bounded (it is the orbit of a bounded set under an action by the unitary group, which is compact.) So the preimage of a bounded set is bounded. It is obvious that the preimage of a closed set is closed.
\end{proof}

We now need some terminology from the theory of real projective curves. Let $\RR\PP^2$ and $\CC\PP^2$ denote the real and complex projective planes respectively. Let $\C \subset  \CC\PP^2$ be a smooth projective curve defined be an equation of degree $n$ with real coefficients. Let $\C(\RR)$ denote $\C \cap \RR\PP^2$. The connected components of $\C(\RR)$ are homeomorphic to circles. It is well known that $\pi_1(\RR\PP^2)$ has two elements. We call a circle embedded in $\RR\PP^2$ an \emph{oval} if it is trivial in $\pi_1(\RR\PP^2)$ and a \emph{pseudoline} if it is nontrivial. An oval $S$ always divides $\RR\PP^2$ into a disk and a M\"{o}bius strip, we call the disk the \emph{interior} of $S$ and the M\"{o}bius strip the \emph{exterior}.

Any two pseudolines in $\RR\PP^2$ intersect, so by our assumption that $\C$ is smooth $\C(\RR)$ contains at most one pseudoline. One can easily see that $\C(\RR)$ contains a pseudoline if and only if $n$ is odd. 

Let $\RR_+\PP^2 \subset \RR\PP^2$ denote the set of $(x:y:z) \in \RR\PP^2$ where $x$, $y$ and $z \geq 0$. 

The main result of Vinnikov is the following:

\begin{VinCrit}
Let $F \in \RR_+[x,y,z]_n$ and let $\C$ be the zero locus of $F$ in $\PP_{\CC}^2$. If $\C$ is smooth, the following are equivalent:

\begin{enumerate}
\item $F$ can be written as $\delta(X,Y,Z)$ with $X$, $Y$ and $Z \in \PD$. In other words, $F \in \HH$. 
\item $\C(\RR)$ consists of precisely $\lfloor n/2 \rfloor$ ovals and, if $n$ is odd, a pseudoline and $\RR_+ \PP^2$ is in the interior of all the ovals.
\item There is a point $(x_0:y_0:z_0) \in \RR_+ \PP^2$ such that every line through $(x_0:y_0:z_0)$ meets $\C$ $n$ times.
\item Every pseudo-line in $\RR\PP^2$ which meets $\RR_+ \PP^2$ intersects $\CC(\RR)$ at least $n$ times.
\end{enumerate}
\end{VinCrit}

\begin{proof}
That (1) and (2) are equivalent is an easy consequence of Theorem 6.1 of \cite{Vin2}. It is topologically clear that (4) is equivalent to (2) and implies (3). By \cite{Vin2}, (3) implies that $F$ can be written as $\det(xX+yY+zZ)$ with $x_0X+y_0Y+z_0Z$ positive definite. However, as $F$ was assumed to have all coefficients positive, it does not vanish on $\RR_+ \PP^2$ and therefore the signature of $xX+yY+zZ$ is constant on $\RR_+ \PP^2$. In particular, $1 \cdot X+0 \cdot Y+0 \cdot Z=X$ is positive definite, as are $Y$ and $Z$.
\end{proof}

\textbf{Remark:}
That (1) implies (4) and hence the other statements is not difficult. A line which meets $\RR_+\PP^2$ can be written as $\{ (a+a't:b+b't:c+c't) \}$ where $(a:b:c)$ and $(a':b':c')$ are two points on that line. We may assume that $(a:b:c)$ and $(a':b':c') \in \RR_+\PP^2$. Then $aX+bY+cZ$ and $a'X+b'Y+c'Z)$ are positive definite; call these matrices $W$ and $W'$. $F$, restricted to the line, is $\det(W+W't)$ and we must show that this polynomial has $n$ real roots. As $W$ is positive definite, we can write $W'=S^* S$, and it is enough to show that $\det((S^{-1})^* W S^{-1} + t \Id)$ has $n$ real roots. The roots of this polynomial are the negatives of the eigenvalues of the positive definite matrix $(S^{-1})^* W S^{-1}$, and hence real.

\begin{Theorem}\label{Torii}
Let $F$ obey the equivalent conditions of the above theorem. In particular, note that the zero locus of $F$ is assumed smooth. Then the fiber of $\delta$ over $F$ is an $\binom{n-1}{2}$ dimensional real torus. 
\end{Theorem}

\begin{proof}
This is more of Theorem~6.1 of \cite{Vin2}.
\end{proof}

\subsection{The Connection Between Horn's Problem and Vinnikov Curves}

In this section, we will describe a new criterion for Horn's multiplicative problem to be solvable in terms of Vinnikov curves.

\begin{prop} \label{NewCrit}
There is an isomorphism (of real semi-algebraic sets) $\beta : \OO^{\mult} \to \P$. Moreover, for $(A,B,C) \in \OO^{\mult}$, if $F=\delta(\beta(A,B,C))$, then the singular values of $A$, $B$ and $C$ are the positive square roots of the zeroes of $F(-1,u,0)$, $F(0,-1,u)$ and $F(u,0,-1)$ respectively.
\end{prop}

\begin{proof}
The map $\beta$ is defined as follows: Choose a representative $(A,B,C)$ for the equivalence class in $\OO^{\mult}$. Choose $P$, $Q$ and $R$ such that $A=PQ^{-1}$, $B=QR^{-1}$ and $C=RP^{-1}$. (One such choice is $(P,Q,R)=(A, \Id, B^{-1})$.) Set $X=P^*P$, $Y=Q^*Q$ and $Z=R^*R$. 

Note that any other choice $(P',Q',R')$ is related to the original choice by $P'=PG$, $Q'=QG$ and $R'=RG$ for some $G \in GL$. Using primes to denote the new variables, $X'=(P')^* P'=G^* P^* P G=G^* X G$, so the map does not depend on the choice of $(P,Q,R)$. 

We must show that it doesn't matter if we choose a different representative  $(A',B',C')=(AU,U^{-1}B,C)$. We again use primes to denote the modified quantities. We may take $P'=P$, $Q'=U^{-1} Q$ and $R'=R$. We have $Y'= (Q')^* Q'=Q^* (U^{-1})^* U^{-1} Q=Q^* Q=Y$. Also, clearly, $X'=X$ and $Z'=Z$.

For any polynomial $f \in \CC(u)$, we denote the roots of $f$ by $\Root(f)$. By definition, $\Crit(A)=\Eigen(A^* A)^{1/2}=\Root(\det (A^* A-u \Id))^{1/2}$. We have
\begin{multline*}
F(-1,u,0)=\det(X-Yu)=\det(P^*P-Q^*Q u)= \\ 
\det((PQ^{-1})^* (P Q^{-1}) - u \Id) \det (Q^* Q)=\det(A^* A-u \Id) \det(Q^* Q).
\end{multline*}

Thus, the roots of $F(-1,u,0)$ are the eigenvalues of $A^* A$.
\end{proof}

The roots of $F(-1,u,0)$, $F(0,-1,u)$ and $F(u,0,-1)$ have geometrical meanings: they correspond to the intersections of $Z(F)$ with the coordinate lines $\{ (x:y:0) \}$, $\{ (0:y:z) \}$ and $\{ (x:0:z) \}$.

We denote the map $\HH \to ((\RR^n_{+})_{\dec})^3$ that takes $F$ to the positive square roots of the zeroes of $F(-1,u,0)$, $F(0,-1,u)$ and $F(u,0,-1)$ by $\partial$. We have just shown:

\begin{NC}
Horn's multiplicative problem is solvable for $(\a,\b,\c)$ if and only if there is an $F \in \HH$ with $\partial F=(\a, \b, \c)$. In other words, if and only if there is a Vinnikov curve passing through $(-1: \alpha_i^2 : 0)$, $(0 : -1 : \beta_i^2)$ and $(\gamma_i^2 : 0 : -1)$.
\end{NC}

\textbf{Remark:} It is obvious that Horn's additive problem is solvable for $(\a, \b, \c)$ if and only if it is solvable for $(2 \a, 2 \b, 2 \c)$. Therefore, using Klyachko's Theorem, Horn's multiplicative problem is solvable for $(\a,\b,\c)$ if and only if it is solvable for $(\a^2, \b^2, \c^2)$. We could have used this observation to eliminate the square roots in the definition of $\partial$. We prefer, however, to reduce our dependence on Klyachko's Theorem as much as possible.

\subsection{Fields of Power Series} \label{Fields}

For technical reasons, we now introduce some fields of power series. Their relevance will become clear in the next section.

Let $\TCC = \bigcup_{n=1}^{\infty} \CC((t^{-1/n}))_{\mathrm{conv}}$ where $\CC((u))$ denotes the field of Laurent series in $u$ and the subscript ``$\mathrm{conv}$'' denotes that we are considering only those series which have a positive radius of convergence around $\infty$. By the leading term of a power series, we mean the term $a t^{\alpha}$ with $\alpha$ the most positive (and $a$ nonzero). Let $\TRR$ denote the corresponding field where the coefficients are real. Let $\TRRp$ be the subsemifield of $\TRR$ where the coefficient of the leading term is positive. It is easy to check that $\TRR$ is an ordered field where the elements of $\TRRp$ are defined to be positive.

\begin{prop}
$\TCC$ is algebraically closed.
\end{prop}

\begin{proof}
This is proven without the hypothesis that the series are convergent in \cite{Wal}, chap. IV, section 3. It is easy but not trivial to modify the proof for the convergent case; a detailed proof can be found in \cite{Pic}, vol. 2, chap. XIII.
\end{proof}

\begin{prop}
$\TRR$ is real closed.
\end{prop}

\begin{proof}
This follows from the Artin-Schrier theorem and the fact that $[\TCC:\TRR]=2$. (See, for example, \cite{Jac}, vol. II, Theorem 11.14.)
\end{proof}

The important property of real closed fields is a result of Tarski's.

\begin{Tar}
A first order statement (a statement constructed from addition, multiplication, equality, inequality and the standard logical connectives and quantifiers) is true in a real closed field if and only if it is true in every real closed field. In particular, a first order statement is true in $\TRR$ if and only if it is true in $\RR$. 
\end{Tar}

\begin{proof}
See, for example, \cite{Jac}, vol. I, sect. 5.6.
\end{proof}

A first order statement is true in $\TRR$ if and only if the first order statement about $\RR$ that results by plugging in $t$ is true for $t$ sufficiently large. (This charcterization wouldn't make sense if we didn't have convergent power series.)

Let $l:\TCC^* \to \QQ$ be the map which sends a power series to the exponent of its leading term. Alternatively, we can define $l(f)$ as $\lim_{t \to \infty} \log f/\log t$. (``$l$'' stands for ``logarithm''.) 

The basic facts about $l$ are the following:

\begin{Lemma}
$l(xy)=l(x)+l(y).$ $l(x+y) \geq \max(l(x), l(y))$, with equality if $l(x) \neq l(y)$ or $x$ and $y \in \TRRp$ (this is not an exhaustive list of cases where equality holds.) If $x$ and $y \in \TRRp$, then $x \leq y$ implies $l(x) \leq l(y)$ and $l(x) < l(y)$ implies $x < y$.
\end{Lemma}

\begin{proof}
Obvious.
\end{proof}

It makes sense to define ``Hermitian'', ``positive definite'', ``eigenvalue'' and so forth for $\TCC$ in exactly the same matter as for $\CC$. By the real-closedness of $\TRR$, all of these definitions will work. (Hermitian matrices will have eigenvalues in $\TRR$, the square roots in the definition of singular values will be defined, etc.) We will denote the power series versions of these concepts by $\GL_n(\TCC)$, $\Her_n(\TCC)$, etc.

The choice of convergent power series with rational exponents was made primarily for convenience at certain technical points. The use of real exponents would require only trivial modifications to our results. The use of formal power series would make it difficult to talk about the toplogies of the curves they defined, which would make it harder to apply Vinnikov's criterion. Our power series are in negative powers of $t$ because this allows us to be consistent with standard sign conventions for hives. 

\subsection{Results}

The philosophy of this paper is the following: $(1/2) \log \HH \in \RR^{\Delta}$ is analogous to $\HIVE$ and we should attempt to make that analogy precise. We do this in two ways: one is by proving statements about power series, the other is by proving statements about approximations. 

The first, easiest part of our analogy is that the maps $\partial : \HH \to (\RR_+^n)_{\dec}^3$ and $\partial_0 : \HIVE \to (\RR^n)_{\dec}^3$ are analogous. We can make this analogy precise in two ways: the first using power series and the second using approximations.

\begin{prop} \label{PowSerBound}
Let $F \in \HH(\TCC)$. Then $l(\partial(F))=(1/2) \partial_0(l(F))$ where $l$ acts on each coordinate.
\end{prop}

\begin{proof}
Let $f(u)=F(u,1,0)$ and write
$$f(y)=\sum f_{n-i} u^i = f_0 \prod (u+r_i^2)$$ 
with $r_1 \geq \cdots \geq r_n$. Unwinding the definitions, we must show that
$$l(r_k)=(1/2) (l(f_{k})-l(f_{k-1})).$$

We have
$$l(f_k)=l \left( \sum_{1 \leq i_1<i_2<\cdots < i_k \leq n} r_{i_1}^2 \cdots r_{i_k}^2 \right).$$
Since every term in the sum is in $\TRRp$, this is just 
$$l( r_1^2 \cdots r_k^2)=2 \sum_{i=1}^k l(r_i).$$
The result is now obvious.
\end{proof}

\begin{prop} \label{RealBound}
Let $F \in \HH(\CC)$. Then $\log(\partial(F))=(1/2) \partial_0(\log F)+O(1)$ where the $O(1)$ depends only on $n$.
\end{prop}

\begin{proof}
This proof is just like the preceeding one except that, instead of the equation
$$l \left( \sum_{1 \leq i_1<i_2<\cdots < i_k \leq n} r_{i_1}^2 \cdots r_{i_k}^2 \right)=l( r_1^2 \cdots r_k^2),$$
we instead notice that
$$\sum_{1 \leq i_1<i_2<\cdots < i_k \leq n} r_{i_1}^2 \cdots r_{i_k}^2 = r_1^2 \cdots r_k^2 C$$
where $1 \leq C \leq \binom{n}{k}$. Thus, 
$$\log \left( \sum_{1 \leq i_1<i_2<\cdots < i_k \leq n} r_{i_1}^2 \cdots r_{i_k}^2 \right) = \log \left( r_1^2 \cdots r_k^2 \right) + O(1).$$

The proof now preceeds as before.
\end{proof}

We now describe the main results of this paper. First we restate the main theorem, which is the result we will actually spend most of our time proving:

\begin{MT}
Let $F = \sum F_{ijk}(t) x^i y^j z^k$ be a homogenous degree $n$ polynomial with coefficients in $\TRRp$ and set $h_{ijk}=l(F_{ijk}) $. If $F$ is a Vinnikov polynomial then $h$ is a hive and, if $h$ is a strict hive then $F$ is a Vinnikov polynomial. 
\end{MT}

From this, we deduce the following, power series free, result.

\begin{Theorem} \label{Real}
There exist vectors $V_1$ and $V_2 \in \RR^{\Delta}$ such that $V_1+\HIVE \supseteq \log \HH \supseteq V_2+\HIVE$.
\end{Theorem}

\begin{proof}
While this statement involves logarithms, it can be restated as a first order statement. The sets $e^{V_r + \HIVE}$ ($r=1,2$) can be described as the subsets of $(\RR_+^n)_{\dec}^3$ where each ratio
$$F_{(i+1)j(k-1)} F_{i(j+1)(k-1)} /  F_{ijk} F_{(i+1)(j+1)(k-2)}$$
is larger than some constant $K_r$ (and similarly for permutations of $i$, $j$ and $k$.) To be explicit, we can restate the result as: ``there exist positive constants $K_1$ and $K_2$ such that, if all of the above ratios are greater than $K_2$, then $F \in \HH$ and, if $F \in \HH$, all of the above ratios are greater than $K_1$.''

We may thus prove this statement in $\TCC$ instead of $\CC$. In this case, the main theorem tells us we may take $K_1=t^{-1}$ and $K_2=t$.
\end{proof}

\textbf{Remark:} Thinking of $\RR_+^{\Delta}=\RR_+[x,y,z]_n$ as the parameter space of degree $n$ plane curves, let $\D \subset \RR_+^{\Delta}$ be the locus corresponding to singular curves. One can check (with a little work) from Vinnikov's criterion that every connected component of $\RR_+^{\Delta} \setminus \D$ either lies completely in or completely out of $\HH$. By proposition \ref{SmoothareDense}, $\HH$ is the closure of a union of connected components of $\RR_+^{\Delta} \setminus D$. The geometry of connected components of $\log \left( \RR_+^{\Delta} \setminus \D \right)$ is studied in chapter 11.5 of \cite{GKZ}. This paper may be thought of as an example of applying the methods of that section. It is not clear how to deduce the results of this paper from those of \cite{GKZ}, although there is a structural similarity.

\begin{Theorem} \label{exists}
If $h_{ijk}$ is a hive then there exists $F=\sum F_{ijk}(t) x^i y^j z^k \in \TRRp[x,y,z]_n$ with $l(F_{ijk})=h_{ijk}$ such that $F \in \HH(\TCC)$.
\end{Theorem}

\begin{proof}
We can take the vectors $V_1$ and $V_2$ in Theorem~\ref{Real} to have coordinates which are logarithms of rational numbers. The statement that those particular vectors have the property of Theorem \ref{Real} is then a first order statement.

Thus, this statement is also true for $\TRRp$. 

Now, suppose that $h_{ijk}$ is a hive and choose $F_{ijk}$ of the form $c_{ijk} e^{t h_{ijk}}$ with $c_{ijk}$. In order to force $F \in \HH(\TCC)$, it is enough to force the ratios 
$$F_{(i+1)j(k-1)} F_{i(j+1)(k-1)} /  F_{ijk} F_{(i+1)(j+1)(k-2)}$$
to be larger than some rational numbers. The fact that $h_{ijk}$ is a hive means that this ratio leads off with a nonnegative power of $t$. If it leads with a positive power of $t$, then it will be greater than any rational number; if it leads off with a constant then, by choosing the $c_{ijk}$ to be the exponential of a point sufficiently inside the hive cone, we can arrange for this ratio to be large enough.
\end{proof}

From the Main Theorem, Theorem~\ref{exists} and Proposition \ref{PowSerBound}, we immediately can immediately deduce:

\begin{Theorem}
Given $(\a, \b, c) \in (\QQ^n_+)^3_{\dec}$, there exist $A$, $B$ and $C \in \GL(\TCC)$ with $l(\Crit(A))=\a$, $l(\Crit(B))=\b$ and $l(\Crit(C))=\c$ if and only if $(\a,\b,\c)$ is the boundary of a hive.
\end{Theorem}

Working a little harder, we can give a new proof of the result of Knutson and Tao:

\begin{HT}
Horn's additive problem is solvable for $(\a, \b, \c)$ if and only if $(\a, \b, \c)$ are the boundary of a hive.
\end{HT}

\begin{proof}
If $S$ and $T$ are two subsets of $\RR^n$, we write $S \sym T$ to denote that there exists a constant $K$ such that every point of $S$ is within $K$ of a point of $T$ and \emph{vice versa}. Note that, if $f$ and $g : \RR^m \to \RR^n$ with $f=g+O(1)$ and $S \subset \RR^m$, then $f(S) \sym g(S)$. Also, if $f$ is linear, then $S \sym T$ implies $f(S) \sym f(T)$.

Theorem~\ref{Real} implies that $\log \HH \sym \HIVE$. As $(1/2) \partial_0$ is linear, $(1/2) \partial_0(\log \HH) \sym (1/2) \partial_0(\HIVE)$. Also, by Proposition~\ref{RealBound}, $\log \circ \partial = (1/2) \partial_0 \circ \log + O(1)$. Thus, 
$$\log \partial(\HH) \sym (1/2) \partial_0(\HIVE)=\partial_0(\HIVE),$$
where the last equality follows simply because $\HIVE$ is a cone.

By the New Criterion, $\partial \HH$ consists of those $(\a, \b, \c)$ for which Horn's multiplicative problem is solvable. By Klyachko's Theorem, $\log \partial(\HH)$ consists of those $(\a, \b, \c)$ for which Horn's additive problem is solvable. 

So the left hand side of the displayed equation consists of those $(\a, \b, \c)$ for which  Horn's additive problem is solvable and the right hand side, by definition, consists of boundaries of hives. However, both are clearly preserved under scaling by $\RR_+$. It is geometrically clear that, if $S$ and $T \subset \RR^n$ are both invariant under $\RR_+$ scaling and $S \sym T$ then $S=T$.
\end{proof}

\textbf{Remark:}
One might try to avoid using Klyachko's Theorem by stating our theorem as ``Horn's multiplicative problem is solvable for $(\a, \b, \c)$ if and only if $(\log \a, \log \b, \log \c)$ is the boundary of a hive.'' But it is not obvious that the set of $(\log \a, \log \b, \log \c)$ such that Horn's multiplicative problem is solvable for $(\a, \b, \c)$ is closed under $\RR_+$ scaling. Thus our proof of the Knutson Tao result uses Klyachko's Theorem in an essential manner, whereas the proof of our other thoerems does not.

\subsection{A Diagram} \label{Diagram}

The following diagram shows most of the spaces and maps discussed in the preceeding section, and maybe useful to refer to.

\xymatrix{
 & \HIVE_n \ar[dr]^{\partial_0} \ar@{^{(}->}[r] & (\RR)^{\Delta_n}/\RR\One  \\
{\OO}^{\add} = \{ (A,B,C) \in (\Her_n \CC)^3 : A+B+C=0 \}/U_n \ar[rr]^-{\Eigen} \ar[d]^{\isomorph}_{\alpha} & & ((\RR^n)_{\dec})^3 \ar[d]^{\mathrm{exp}}_{\isomorph} \\
{\OO}^{\mult} = \{ (A,B,C) \in (\GL_n \CC)^3 : ABC=1 \}/U_n^3 \ar[rr]^-{\Crit} \ar[d]^{\isomorph}_{\beta} & & ((\RR_+^n)_{\mathrm{dec}})^3 \ar@{=}[d]  \\
{\P} = \{ (X,Y,Z) \in (\PD_n \CC)^3 \}/ \GL_n \CC \ar@{>>}[r]^-{\delta}_-{\det(xX+yY+zZ)} & {\HH}_n \ar[r]^-{\partial}\ar@{^{(}->}[d]  & ((\RR_+^n)_{\mathrm{dec}})^3 \\
& \RR_+[x,y,z]_n/\RR_+ \ar@{=}[r] &  (\RR_+)^{\Delta_n}/\RR_+ 
}

\section{Patchworking}

Patchworking is a technique developed by Viro to study the topology of families of real plane curves $\sum f_{ijk}(t) x^i y^j z^k=0$ for $t$ large. This section is dedicated to describing the results we will need from this theory.

\subsection{Viro's Theorem for Triangulations}

Let $\alpha_{ijk} \in \RR^{\Delta}$ and $f_{ijk} \in \TRR$ with $l(f_{ijk})=\alpha_{ijk}$.  Set 
$$F(t)(x,y,z)=\sum f_{ijk} x^i y^j z^k.$$
Let $\C_t$ denote the zero locus of $F$ evaluated at $t \in \RR_{+}$. (This is defined for $t$ sufficiently large.)

The function $\alpha_{ijk}$ gives rise to a polyhedral subdivision of $\Delta$ according to the following prescription: consider the convex hull of the set of points $(i,j,k,\alpha_{ijk}) \in \Delta \times \RR \subset (\RR^3/\RR) \times \RR$. Take the ``upper'' faces of this convex hull and project them down to $\Delta$. (For more on this procedure, see chapter 7 of \cite{GKZ}.) If $\alpha$ is chosen generically, this subdivision will be a triangulation.

The following is obvious.
\begin{prop}\label{HivesandTriangulations}
$\alpha$ is a strict hive if and only if the corresponding triangulation is the standard triangulation of $\Delta$ (into $n^2$ equilateral triangles of side length 1.) $\alpha$ is a hive if and only if the corresponding polyhedral subdivision is a coarsening of the standard subdivision. 
\end{prop}

\begin{figure}
\centerline{\includegraphics{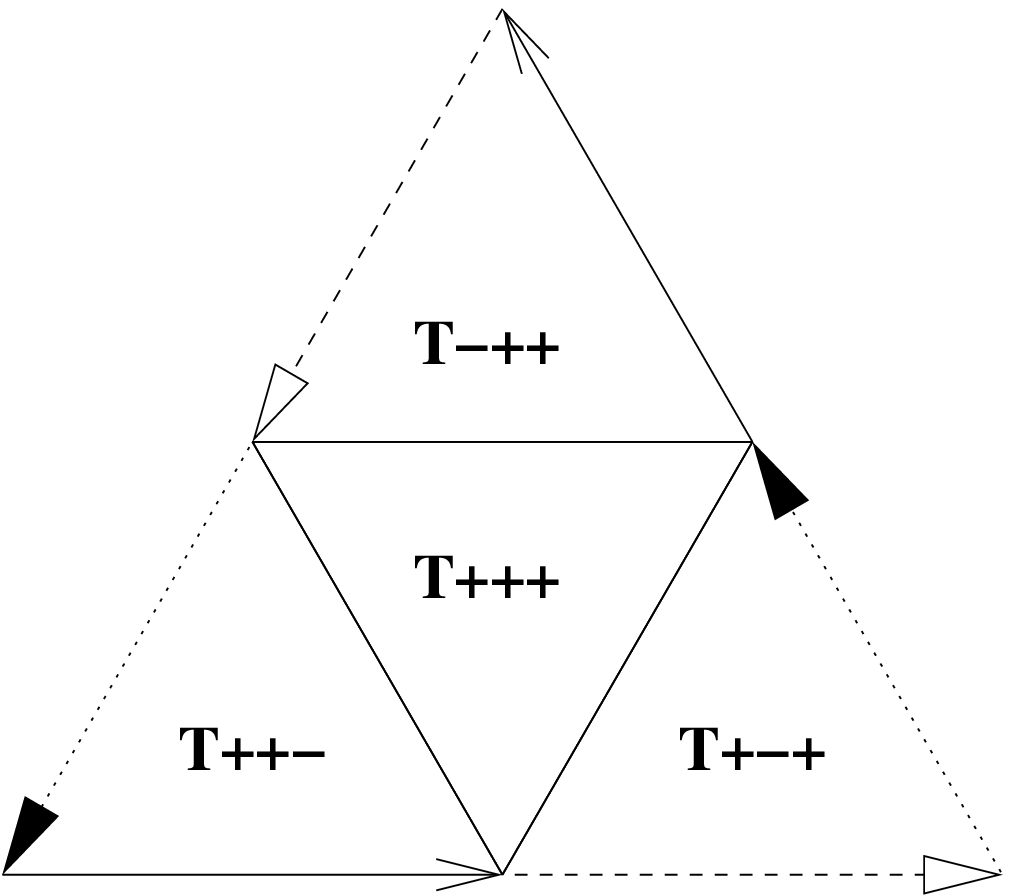}}
\caption{How to Glue the $\Tau$'s into $\RR\PP^2$.}\label{Gluing}
\end{figure}

Let $(\epsilon_1, \epsilon_2, \epsilon_3) \in \{ 1, -1 \}^3$. Define $\Tau_{\epsilon_1, \epsilon_2, \epsilon_3} \subset \RR\PP^2$ to be the set of $(x_1:x_2:x_3)$ such that $\epsilon_i x_i \geq 0$. We will usually abbreviate $\Tau_{1, 1, 1}$ by $\Tau_{+++}$ and so forth. Note that $\Tau_{+++}$ is what we previously called $\RR_+ \PP^2$ and $\Tau_{\epsilon_1, \epsilon_2, \epsilon_3}=\Tau_{-\epsilon_1, -\epsilon_2, -\epsilon_3}$. Each $\Tau_{\epsilon_1, \epsilon_2, \epsilon_3}$ is a topological disk and should be visualized as a triangle. Our goal will be to describe the topology of $\C \cap \Tau_{\epsilon_1, \epsilon_2, \epsilon_3}$. Pasting together the answer for each $(\epsilon_1, \epsilon_2, \epsilon_3)$, we will find the topology of $\C \cap \RR\PP^2=\C(\RR)$. The gluing of the $\Tau$'s is shown in Figure~\ref{Gluing}.

Fix $(\epsilon_1, \epsilon_2, \epsilon_3) \in \{ 1, -1 \}^3$. Suppose that $\alpha$ induces a triangulation $T$ of $\Delta$.  Label each vertex $(i,j,k) \in \Delta$ with either a $+$ or a $-$ by giving $(i,j,k)$ the sign of
$\epsilon_1^i \epsilon_2^j \epsilon_3^k f_{ijk}.$

We will now describe a collection of polygonal paths in the convex hull of $\Delta$: if $((i,j,k), (i',j',k'),(i'',j'',k''))$ is a triangle of $T$ such that $(i,j,k)$ and $(i',j',k')$ have the same sign and $(i'',j'',k'')$ has the opposite sign, draw a line segment connecting the midpoints of edges $((i,j,k),(i'',j'',k''))$ and $((i',j',k'), (i'',j'',k''))$. If all three vertices have the same sign, do not draw any segments.

\begin{VT}
For $t$ large enough, there is a homeomorphism between $\Tau_{\epsilon_1, \epsilon_2, \epsilon_3}$ and the convex hull of $\Delta$, taking the corresponding sides of the two triangles to each other and taking the polygonal paths described above to $\C \cap \Tau_{\epsilon_1, \epsilon_2, \epsilon_3}$.
\end{VT}

\begin{proof}
This is essentially proven in \cite{Viro}. (Viro uses monomials rather than arbitrary members of $\TRR$, but one can check that this is unimportant.)
\end{proof}

\subsection{Viro's Theorem for Arbitrary Polyhedral Subdivisions}

Now suppose that $\alpha$ induces a polyhedral subdivision $\Delta=\bigcup_{\Gamma \in \Phi} \Gamma$ which is not necessarily a triangulation. Let $\Phi$ denote the set of faces of the subdivision; we emphasize that $\Phi$ contains the edges and vertices of the subdivision and not only the two dimensional faces. Note that not all the vertices of $\Delta$ need be used in $\Phi$.

Let $\Gamma \in \Phi$. There is a linear functional $\alpha^{\Gamma} : \RR^3 / \RR \to \RR$ such that $\alpha_{\Gamma}(i,j,k)=\alpha(i,j,k)$ if and only if $(i,j,k) \in \Gamma$. (If $\Gamma$ is not two dimensional, $\alpha_{\Gamma}$ is not uniquely defined, but this will not be important.) 

For $(i,j,k) \in \Delta$, write 
$$f_{ijk}=c^{\Gamma}_{ijk} t^{\alpha^{\Gamma}(i,j,k)} (1+o(1))$$
where $c_{ijk} \in \RR$ and $o(1)$ denotes a member of $\TRR$ with $l(o(1)) < 0$. We have $c_{ijk} \neq 0$ if and only if $(i,j,k) \in \Gamma$. Set $F^{\Gamma}=\sum c_{ijk}^{\Gamma} x^i y^j z^k$. Intuitively, one should think of $F^{\Gamma} (x,y,z)$ as a good approximation to $F(t^{-\alpha^{\Gamma}(1,0,0)} x, t^{-\alpha^{\Gamma}(0,1,0)} y, t^{-\alpha^{\Gamma}(0,0,1)} z)$ for $(x,y,z)$ fixed and $t$ large. 

We can consider $F^{\Gamma}$ as a function on $(\CC^*)^{\dim \Gamma}$. (More properly, we should consider $F^{\Gamma}$ as a section of a line bundle on the toric variety corresponding to $\Gamma$, but this would introduce concepts we have no need for.) $\Gamma$ and $\RR_+^{\dim \Gamma}$ are clearly both topologically disks of dimension $\dim \Gamma$.  There is a correct way to choose such homeomorphisms so that the curves $Z(F^{\Gamma}) \cap \RR_+^{\dim \Gamma}$ paste together to form a closed subspace of $\Delta$.

The following result follows from Viro's proof:

\begin{Theorem}
If all of the curves $Z(F^{\Gamma}) \cap (\CC^*)^{\dim \Gamma}$ are smooth then, for $t$ large enough, there is a homeomorphism from $\Tau_{+++}$ to the convex hull of $\Delta$ taking $Z(F) \cap \Tau_{+++}$ to the subspace of $\bigcup_{\Gamma \in \Phi} \Gamma$ constructed above. Similar comments apply to the other $\Tau_{\epsilon_1, \epsilon_2, \epsilon_3}$'s, with appropriate changes of sign in the above construction.
\end{Theorem}

It is clear that Viro's patchworking theorem is a special case of the above.

The case where the curves $Z(F^{\Gamma})$ are not smooth is still more complicated to describe. In this case, one can not completely describe the topology of $Z(F)$ in terms of the $Z(F^{\Gamma})$ but one can come close. Namely, one can find a homeomorphism $\iota$ as in the above theorem such that, for $t$ large enough, the topology of $\iota(Z(F)) \cap \Gamma$ is a perturbation of the topology of $Z(F^{\Gamma}) \cap \RR_+^{\dim \Gamma}$.

We mention this only for the special case where $\Gamma$ is one dimensional which will be important later.

\begin{Lemma} \label{SignChange}
With the notation above, suppose that there is a path $\gamma$ through $\Delta$ along the edges (one dimensional faces) of $\Phi$ such that the sign of $(i,j,k)$ changes $m$ times along $\gamma$. Then there is a homeomorphism $\Delta \to \Tau_{+++}$ such that, under this isomorphism, $\gamma$ meets $Z(F)$ at most $m$ times. Similar results hold for the other $\Tau_{\epsilon_1 \epsilon_2 \epsilon_3}$'s.
\end{Lemma}

\begin{proof}
After the preceeding discussion, it is enough to show that, if there are only $m$ sign changes in the coefficients of $f(u) \in \RR[u]$ then $f$ has at most $m$ positive roots. This is Descartes' rule of signs.
\end{proof}

\section{Proof of the Main Theorem}

The aim of this section is to prove the Main Theorem, which we restate here for the reader's convenience.

\begin{MT}
Let $F=\sum F_{ijk}(t) x^i y^j z^k \in \TRRp[x,y,z]_n$. Set $l(F_{ijk})=h_{ijk}$. If $h_{ijk}$ is a strict hive then $F \in \HH(\TCC)$. If $F \in \HH(\TCC)$ then $h_{ijk}$ is a hive.
\end{MT}

In both cases, the key to the proofs will be to combine Vinnikov's Criterion, Viro's patchworking methods and a little combinatorial reasoning about triangulations.

\subsection{If $h_{ijk}$ is a Strict Hive then $F \in \HH(\TCC)$}

By Proposition~\ref{HivesandTriangulations}, $h$ induces the standard subdivision of $\Delta$ and by assumption all of the coefficients of $F$ are positive. By the Vinnikov Criterion, it is enough to show that $Z(F(t)) \cap \RR\PP^2$ has the correct topology for $t$ large enough. Thus, we simply must carry out the construction in Viro's Theorem for the standard triangulation with all of the $F_{ijk}$ positive.

\begin{figure}
\centerline{\includegraphics{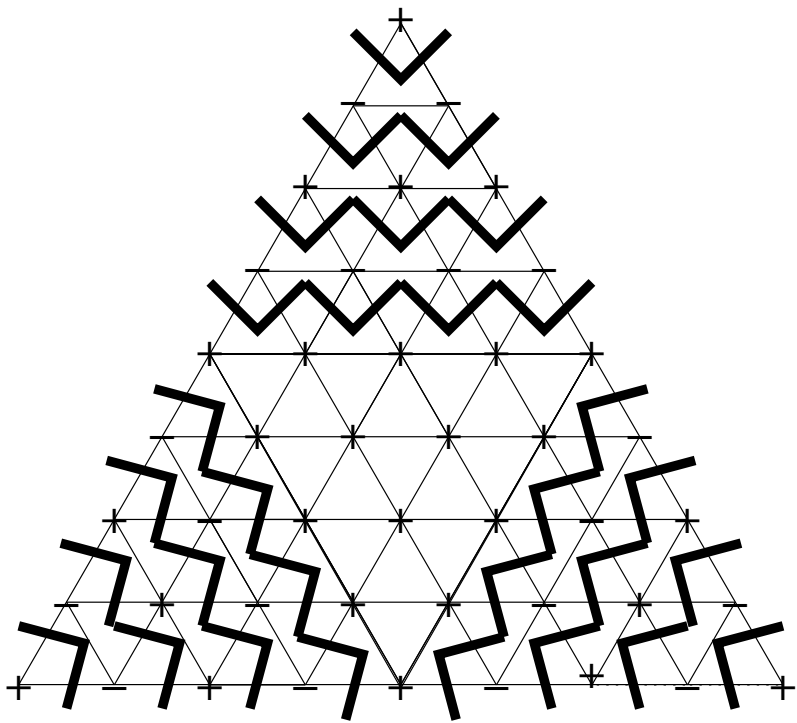}}
\caption{The Standard Triangulation of $\Delta$ and the Topology of the Resulting Curve} \label{Standard}
\end{figure}

The result is shown in Figure~\ref{Standard}, where the boundary of the triangle is glued to itself to form $\RR\PP^2$ as shown in Figure~\ref{Gluing} and the bold lines indicate the polygonal paths.

\subsection{If $F \in \HH(\TCC)$ then $h_{ijk}$ is a Hive}

Suppose for the sake of contradiction that $h_{ijk}$ is not a hive. By lemma~\ref{SmoothareDense} we may perturb each $F_{ijk}$ such that $F$ is a smooth curve and still in $\HH$. If we perturb $F_{ijk}$ by less than $t^{h_{ijk}}$, this will not change the $h$'s. Thus, without loss of generality, we may assume that $F$ is smooth (over $\TCC$ or, equivalently, that $F(t)$ is smooth for $t$ sufficiently large. So Vinnikov's criterion applies to $F(t)$ for $t$ sufficiently large.

 Let $\Phi$ be the polyhedral subdivision induced by $h$. Since $\Phi$ is not a coarsening of the standard triangulation, it contains an edge $e$ which is not in the standard triangulation. Let this edge run between $(i_0,j_0,k_0)$ and $(i_1,j_1,k_1)$. At least one of $|i_0-i_1|$, $|j_0-j_1|$ and $|k_0-k_1|$ must be greater than $1$. Without loss of generality, suppose that $i_0-i_1 > 1$.

\begin{Lemma}
For some $j_2$ and $k_2$ with $j_2+k_2=n$, there is a path from $(n,0,0)$ to $(0,j_2,k_2)$ traveling along the edges of $\Phi$ and using no more than $n-1$ edges.\end{Lemma}

\begin{proof}
By starting at $(i_0,j_0,k_0)$ and traveling along the edges of $T$ in the direction of increasing $i$, we may find a path from $(i_0,j_0,k_0)$ to $(n,0,0)$ which uses no more than $n-i_0$ edges of $\Phi$. Similarly, traveling in the direction of decreasing $i$, we may find a path from $(i_1,j_1,k_1)$ to $(0,j_2,k_2)$ which uses no more than $i_1$ edges of $\Phi$. Concatenating these two paths and the edge $e$, we have a path of length no more than $(n-i_0)+1+i_1 \leq n-1$.
\end{proof}

\begin{figure}
\centerline{\includegraphics{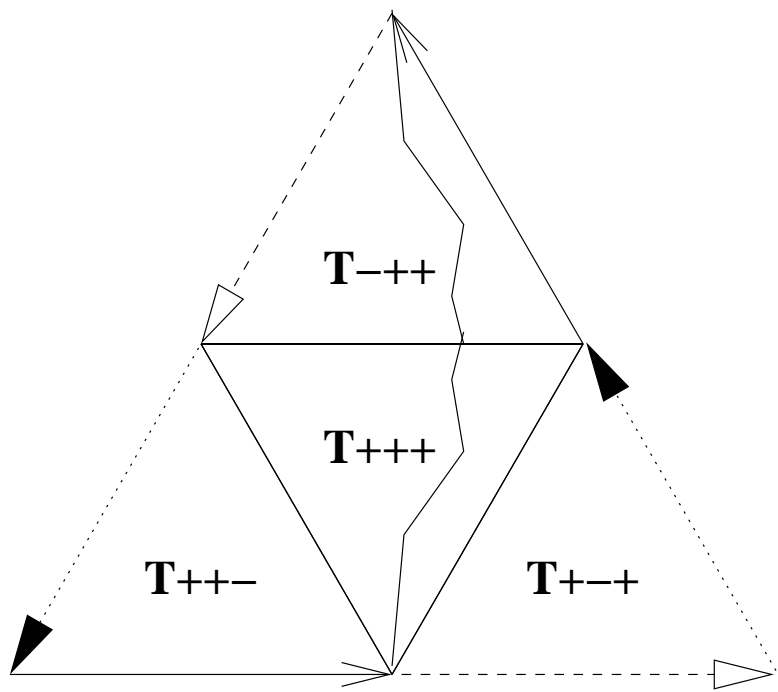}}
\caption{The Pseudo-Line Made Up of Two Copies of $\gamma$.} \label{Loop}
\end{figure}

Let $\gamma$ be the path guaranteed by the lemma. Build a model of $\RR\PP^2$ by gluing together the $\Tau$'s and identify each $\Tau$ with $\Delta$ as in Viro's Theorem. Then the concatenation of the images of $\gamma$ in $\Tau_{-++}$ and $\Tau_{+++}$ form a pseudoline. See Figure~\ref{Loop} for an illustration. 
The former copy of $\gamma$ has at most $n-1$ sign changes, as it has at most $n-1$ edges. The latter copy lies in $\Tau_{+++}$, where all signs are the same, and thus has no sign changes. So this pseudoline has at most $n-1$ sign changes along it. From Lemma~\ref{SignChange}, for $t$ sufficiently large, there is a pseudoline in $\RR\PP^2$, meeting $\RR_+ \PP^2$ and crossing $Z(F)$ at most $n-1$ times. This contradicts Vinnikov's criterion for $F$ to lie in $\HH$. 

\section{An Explicit Value for $V_1$}

We have seen that the Main Theorem combined with a little logic implies that there exist vectors $V_1$ and $V_2 \in \RR^{\Delta}$ such that $V_1+\HIVE \supseteq \log \HH \supseteq V_2+\HIVE$. In this section, we will find an explicit value for $V_1$. Our proof will rely on Vinnikov's criterion but not on patchworking or the Main Theorem, so it may be used to give an alternative proof of the second direction of the Main Theorem.

Specifically, we will prove that:

\begin{prop} \label{Backward}
Let $\sum F_{ijk} x^i y^j z^k \in \HH$. Then

\begin{eqnarray*}
\frac{ 2 (k-1) }{k} F_{(i+1)j(k-1)} F_{i(j+1)(k-1)} & > & F_{ijk} F_{(i+1)(j+1)(k-2)} \\
\frac{ 2 (j-1) }{j} F_{(i+1)(j-1)k} F_{i(j-1)(k+1)} & > & F_{ijk} F_{(i+1)(j-2)(k+1)} \\
\frac{ 2 (i-1) }{i} F_{(i-1)(j+1)k} F_{(i-1)j(k+1)} & > & F_{ijk} F_{(i-2)(j+1)(k+1)}
\end{eqnarray*}
\end{prop}

We have the following corollary:

\begin{cor}
In the above notation, one may take 
$$V_1(i,j,k)= - \log \left( i! j! k! 2^{ij+jk+ki} \right)$$
\end{cor}

\begin{proof}
Suppose that $\sum F_{ijk} x^i y^j z^k \in \HH$. We must show that $\log F_{ijk} - V_1(i,j,k)$ is a hive.

With the above $V_1$, we have
$$ \frac{e^{V_1(i,j,k)} e^{V_1(i+1,j+1,k-2)}}{e^{V_1(i,j+1,k-1)} e^{V_1(i+1,j,k-1)}} = \frac{2 (k-1)}{k}.$$

Thus, the first line of the conclusion of Proposition~\ref{Backward} can be restated as
$$\left( \frac{ F_{(i+1)j(k-1)} }{e^{V_1(i+1,j,k-1)}} \right) \left( \frac{ F_{i(j+1)(k-1)} }{e^{V_1(i+1,j,k-1)}} \right) >
\left( \frac{ F_{ijk} }{e^{V_1(i,j,k)}} \right) \left( \frac{ F_{(i+1)(j+1)(k-2)} }{e^{V_1(i+1,j+1,k-2)}} \right)$$
and similarly for the other lines. Taking $\log$s of both sides, we conclude that $\log F_{ijk} - V1(i,j,k)$ is indeed a hive.
\end{proof}

\subsection{Derivatives of Vinnikov Curves}

The aim of this section is to prove the following technical lemma:

\begin{Lemma}\label{Deriv}
Let $F(x,y,z) \in \HH_n$. (The subscript $n$ means that we are dealing with degree $n$ polynomials.) Let $x_0$, $y_0$ and $z_0 \geq 0$, $(x_0,y_0,z_0) \neq (0,0,0)$. Then 
$$x_0 \frac{\partial F}{\partial x} + y_0 \frac{\partial F}{\partial y} + z_0 \frac{\partial F}{\partial z} \in \HH_{n-1} .$$
\end{Lemma}

\begin{proof}
By lemmas~\ref{SmoothareDense} and~\ref{proper} we may assume that $F$ is smooth.

Fix $(a:b:c) \in \RR\PP^2 \setminus \{ (x_0:y_0:z_0) \}$. For $t \in \RR$, set $l(t)=(a+x_0t, b+y_0t,c+z_0t)$. So $l(t)$ traces out the line in $\RR\PP^2$ joining $(x_0:y_0:z_0)$ and $(a:b:c)$. By Vinnikov's criterion, $l(t)$ meets $Z(F)$ at $n$ real points. In other words, the polynomial $f(t)=F(a+x_0t,b+y_0t,c+z_0t)$ has $n$ distinct real roots. By Rolle's theorem, $df/dt$ has $n-1$ distinct real roots.

But an easy computation shows that 
$$ \left. \frac{df}{dt} \right|_{t=t_0} = \left. \left( x_0 \frac{\partial F}{\partial x} + y_0 \frac{\partial F}{\partial y} + z_0 \frac{\partial F}{\partial z} \right) \right|_{(x,y,z)=l(t_0)}.$$
So we have shown that $x_0 \partial F/\partial x+y_0 \partial F/\partial y+z_0 \partial F/\partial z$ has $n-1$ distinct real zeroes on the line through $(a:b:c)$ and $(x_0:y_0:z_0)$.

But $(a:b:c)$ was chosen arbitrarily. So we have shown that every line through $(x_0:y_0:z_0)$ meets $Z(x_0 \partial F/\partial x+y_0 \partial F/\partial y+z_0 \partial F/\partial z)$ at $n-1$ distinct points. Also, it is clear that $x_0 \partial F/\partial x+y_0 \partial F/\partial y+z_0 \partial F/\partial z$ has positive coefficients. So, by Vinnikov's criterion, $x_0 \partial F/\partial x+y_0 \partial F/\partial y+z_0 \partial F/\partial z \in \HH_{n-1}$.
\end{proof}

\subsection{The Proof}

We will now prove Theorem~\ref{Backward}. We will show 
$$\frac{ 2 (k-1) }{k} F_{(i+1)j(k-1)} F_{i(j+1)(k-1)} > F_{ijk} F_{(i+1)(j+1)(k-2)},$$
the proofs of the other inequalities are similar.

Our proof is by induction on $n$. 

Our base case is $n=2$, so $(i,j,k)=(0,0,2)$. We may assume without loss of generality that $X=\Id$, put $Y=\left( \begin{smallmatrix} Y_{11} & Y_{12} \\ Y_{21} & Y_{22} \end{smallmatrix} \right)$ and $Z=\left( \begin{smallmatrix} Z_{11} & Z_{12} \\ Z_{21} & Z_{22} \end{smallmatrix} \right)$. We are being asked to show that
$$(Y_{11}+Y_{22}) (Z_{11}+Z_{22}) \geq Y_{11} Z_{22}+Y_{22} Z_{11}$$
or, in other words, that
$$Y_{11} Z_{11}+Y_{22} Z_{22} > 0$$
As the diagonal terms of $Y$ and $Z$ are positive, this is obvious.

We now continue with the induction. Suppose that $i+j+k=n \geq 3$. Our proof divides into three cases:

\textbf{Case 1:} $i > 0$.

Set $F'=\partial F/\partial x$, so $F'(x,y,z)=\sum i F_{ijk} x^{i-1} y^j z^k$. By lemma~\ref{Deriv}, $F' \in \HH_{n-1}$. By induction,
$$\frac{ 2 (k-1) }{k} \left( (i+1) F_{(i+1)j(k-1)} \right) \left( i F_{i(j+1)(k-1)} \right) > \left( i F_{ijk} \right) \left( (i+1) F_{(i+1)(j+1)(k-2)} \right).$$
Cancelling $i(i+1)$ from both sides, we are done.

\textbf{Case 2:} $j>0$.

This case is precisely analogous to the previous case.

\textbf{Case 3:} $k > 2$.

Set $F'=\partial F/\partial z$, so $F'(x,y,z)=\sum k F_{ijk} x^i y^j z^{k-1}$. By lemma~\ref{Deriv}, $F' \in \HH_{n-1}$. By induction,
$$\frac{ 2 (k-2) }{k-1} \left( (k-1) F_{(i+1)j(k-1)} \right) \left( (k-1) F_{i(j+1)(k-1)} \right) > \left( k F_{ijk} \right) \left( (k-2) F_{(i+1)(j+1)(k-2)} \right).$$
Dividing $k(k-2)$ out of each side, we are done. \qedsymbol

\section{Connection to Honeycombs and Amoebae} \label{Honey}

\subsection{Amoebae}

Let $F \in \CC[x,y,z]$ be a polynomial. There is a map $\log : (\CC^*)^3 \to \RR^3$ by $(z_1, z_2, z_3) \mapsto (\log |z_1|, \log |z_2|, \log |z_3|)$. The \emph{amoeba} of $F$, which we will denote by $\A(F)$, is defined to be $\log \left( Z(F) \cap (\CC^*)^3 \right)$. If $F$ is homogenous, $\A(F)$ is preserved under translation by $(1,1,1)$ and we will abuse notation by using the same sybol to refer to the image of $\A(F)$ in $\RR^3/(1,1,1)$. We will not aim to discuss the theory of amoebae deeply, see \cite{Mik} for more background.

\subsection{NonArchimedean Amoebae}

Similarly, let $F \in \TCC[x,y,z]$. The \emph{tropicalization}, also known as the \emph{nonarchimedean amoeba}, of $F$ is $l \left( Z(F) \cap (\TCC^*)^3 \right)$. We denote the tropicalization of $F$ by $\Trop F$. 

There is a second description of $\Trop F$. Write $F=\sum f_{ijk} x^i y^j z^k$.  Define a piecewise linear functions $\phi: \RR^3 \to \RR$ by $\phi(x,y,z)=\max_{(i,j,k) \in \Delta} (v(f_{ijk})+xi+yj+zk)$. $\Trop F$ is the nonsmooth locus of this function.

The connection between nonarchimedean amoebae and ordinary amoebae is the following:

\begin{prop}
Let $F \in \TCC[x,y,z]$. 
$$\lim_{t \to \infty} \frac{1}{\log t} \A(F(t)) = \Trop F$$
where the convergence is in the Hausdorff metric. 
\end{prop}

\subsection{Honeycombs}

Now, suppose that $F \in \TCC[x,y,z]_n$ is a Vinnikov curve, so that $v(f_{ijk})$ form a hive. Then $\Trop F$ is what is known as the \emph{honeycomb of $h$}. 

\begin{figure}
\centerline{\includegraphics{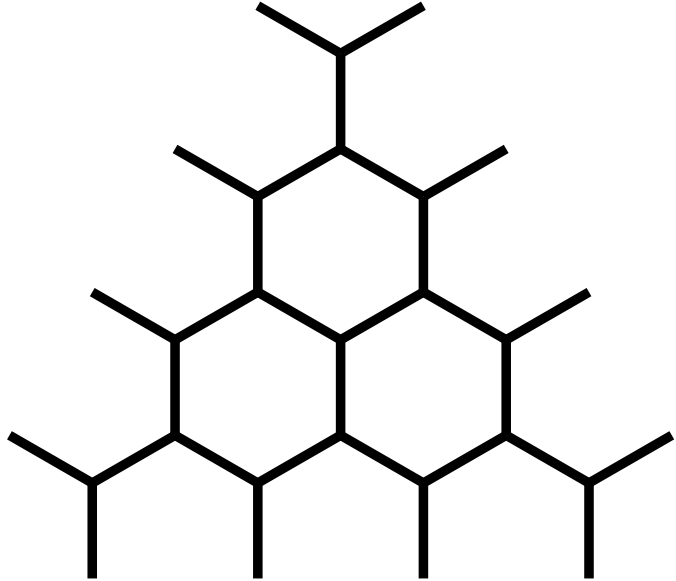}}
\caption{A Honeycomb Arising from a Hive of Order $4$.} \label{honey}
\end{figure}

If $h$ is a strict hive, its honeycomb will consist of a grid of hexgons, as shown in figure~\ref{honey}. If $h$ is a hive which is not strict, some of the edge lengths in the honeycomb will degenerate to $0$; see \cite{KT} for a precise statemement of the sort of degernerations that can occur. Giving a honeycomb is precisely equivalent to giving a hive modulo $\One$. 

The boundary of $h$ is simply the values of $x-y$, $y-z$ and $z-x$ on the unbounded rays of $H(h)$. So, the Hive Theorem says that Horn's additive problem is solvable for $\a$, $\b$ and $c$ if and only if there is a honeycomb whose unbounded rays are at positions corresponding to $\a$, $\b$ and $\c$. 

\begin{figure}
\centerline{\includegraphics{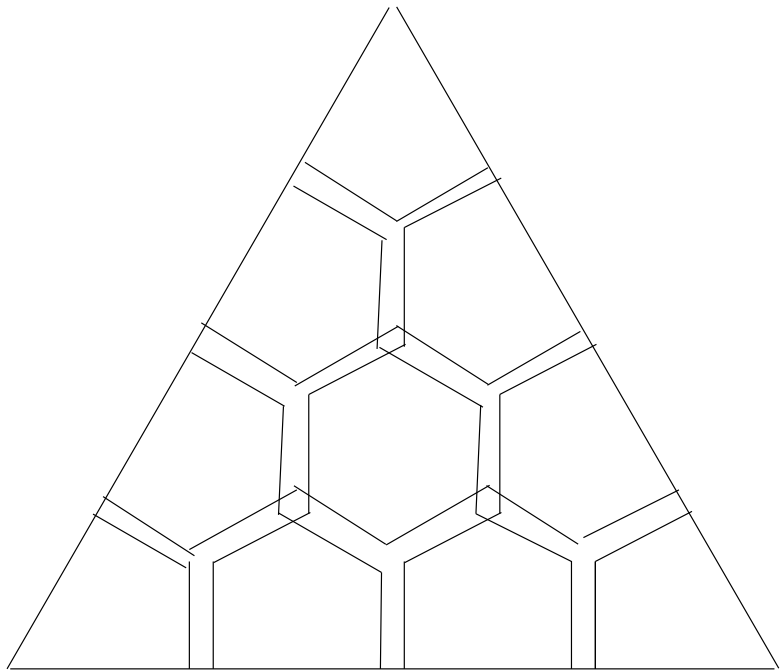}}
\caption{The Real Ovals of a Vinnikov Curve Superimpose to form a Honeycomb.} \label{fold}
\end{figure}

It is interesting to understand how the hexagonal geometry of the honeycomb relates to the nested curve geometry of the Vinnikov curve. The answer is that the boundary of the amoeba of a Vinnikov curve is made up almost entirely of the real points of the Vinnikov curve except for some segments of length $o(\log t)$.  In the limit as $t \to \infty$, the limit of the real points of $Z(F)$ is the same as the limit of the whole amoeba.  Thus, $\Trop F$ is the limit of the image of the real points of $Z(F)$ under the $\log |\cdot|$ map. This map folds the jagged arcs of the Vinnikov curve in figure \ref{Standard} up into a honeycomb, as shown in figure \ref{fold}. In figure \ref{fold}, the lines should be precisely superimposed to yield a honeycomb, but have been seperated slightly for clarity.

I am grateful to Mikhalkin for helping me understand the appearance of these amoebae more closely.

\section{Behavior in Direct Sums} \label{Sum}

Let $\psi$ denote the composite map $\delta \circ \beta : \OO^{\mult}_n \to \RR_+[x,y,z]_n$. Let $(A,B,C)$ and $(A',B',C')$ lie in $\OO^{\mult}_n$ and $\OO^{\mult}_{n'}$ respectively. Then we can consider $(A \oplus A', B \oplus B', C \oplus C') \in \OO^{\mult}_{n+n'}$. Then the following is obvious

\begin{prop}
$$\psi(A,B,C) \psi(A',B',C') = \psi (A \oplus A', B \oplus B', C \oplus C')$$
\end{prop}

The map $\alpha : \OO^{\add} \to \OO^{\mult}$ constructed in \cite{AMW} preserves direct sums so this also holds for the map $\OO^{\add} \to \RR_+[x,y,z]_n$.

Let $H$ and $H'$ lie in $\HIVE_n$ and $\HIVE_{n'}$ respectively. We define the convolution $H * H' \in \HIVE_{n+n'}$ to be the hive defined by
$$(H * H')_{I,J,K}=\max_{\substack{i+i'=I \\ j+j'=J \\ k+k'=K}} (H_{ijk} + H'_{i'j'k'}).$$
This was proved to be a hive in, for example, $\cite{DK}$ and corresponds to the overlay of honeycombs. It is widely accepted that overlay of honeycombs should correspond to direct sum of matrices.

We now check that the above proposition is consistent with the $*$ operator. We state the comparison in terms of power series; there is no difficulty in proving a similar asymptopic result without power series.

\begin{prop}
Let $(A,B,C)$ and $(A',B',C')$ lie in $\OO^{\mult}_n(\TCC)$ and $\OO^{\mult}_{n'}(\TCC)$. Then
$$l(\psi(A \oplus A', B \oplus B', C \oplus C')) = l(\psi(A,B,C)) * l(\psi(A',B',C')).$$
\end{prop}

\begin{proof}
Letting $F$, $f$ and $f'$ denote $\psi(A \oplus A',B \oplus B',C \oplus C')$, $\psi(A,B,C)$ and $\psi(A',B',C')$ respectively, we observed before that $F=ff'$.  Thus, we have
$$F_{IJK}=\sum_{\substack{i+i'=I \\ j+j'=J \\ k+k'=K}} f_{ijk} f'_{i'j'k'}$$
and
$$l(F_{IJK})=l \left( \max_{\substack{i+i'=I \\ j+j'=J \\ k+k'=K}} (f_{ijk} f'_{i'j'k'}) \right)$$
(Here we used that the terms in the sum are all positive.) So
$$\log(F_{IJK})=\max_{\substack{i+i'=I \\ j+j'=J \\ k+k'=K}} \left( \log(f_{ijk}) + \log(f'_{i'j'k'}) \right)$$
\end{proof}

\section{Future Directions}

\subsection{The Ronkin Function}

In this section, we will describe a map $u:\RR[x,y,z]_n \to \RR^{\Delta}$ which has proven to be of use in the theory of amoebas. We will describe reasons to believe that it might be profitable to modify the diagram in section~\ref{Diagram} by replacing the coordinate-wise $\log$ by $u$. For more background, see $\cite{Mik}$. Kenyon and Okounkov (see \cite{KO}) have recently had great success in parameterizing the space of Hanack curves, another type of plane curve with specified topology, by using these methods.

Let $F \in \RR[x,y,z]_n$ (or $\CC[x,y,z]_n$). Let $(x,y,z) \in \RR^3$. We define
$$N_F(x,y,z)= \frac{1}{(2 \pi)^3} \int_0^{2 \pi} \int_0^{2 \pi} \int_0^{2 \pi} \log \left| F(e^{x+i \theta_1}, e^{y+i \theta_2}, e^{z+i \theta_3}) \right| d\theta_1 d\theta_2 d\theta_3.$$
Let $(i,j,k) \in \Delta$, we put
$$u_{ijk}(F)=\inf_{(x,y,z) \in \RR^3} \left( N_F(x,y,z)-ix-jy-kz \right).$$

The integral defining $N_F$ always converges, even when the quantity inside the logarithm is sometimes 0. $N_F$ is a convex function. The infimum in the definition of $u_{ijk}$ is always finite and attained as long as $F_{n00}$, $F_{0n0}$ and $F_{00n}$ are nonzero. 

The following result shows that $u_{ijk}$ and $\log F_{ijk}$ have the same asymptotic behavior in the case that interests us. 

\begin{prop}
Let $F(t)(x,y,z)=\sum F_{ijk} x^i y^j z^k$ with $F_{ijk} \in \TCC$, $l(F_{ijk})=\alpha_{ijk}$. Let $\Phi$ be the subdivision of $\Delta$ induced by $\alpha_{ijk}$. Suppose that $(i,j,k)$ is a vertex of a facet of $\Phi$. Then 
$$\lim_{t \to \infty} u_{ijk}(F(t))/\log t = \lim_{t \to \infty} \log |F_{ijk}(t)|/\log t=\alpha_{ijk}.$$

\end{prop}

\begin{proof}
This is an easy consequence of Theorem 2 of \cite{PR}.
\end{proof}

Note that the hypothesis of the theorem is satisfied for every $(i,j,k) \in \Delta$ when $\alpha_{ijk}$ is a strict hive. We now prove two Propositions that suggest $u_{ijk}$ may be more important that $\log F_{ijk}$. We let $\psi' : \OO^{\mult} \to \RR^{\Delta}$ be the map defined by $u \circ \delta \circ \beta$.

\begin{prop}\label{SuperPosRankin}
Let $(A,B,C) \in \OO^{\mult}_n$ and $(A',B',C') \in \OO^{\mult}_{n'}$. Then 
$$\psi'(A \oplus A', B \oplus B', C \oplus C') \geq \psi'(A,B,C) * \psi'(A',B',C').$$
(The notation $*$ is defined in section \ref{Sum}.) 

If $(A,B,C)$ and $(A',B',C')$ lie in $\OO^{\mult}(\TCC)$ instead, then, moreover, we have equality for sufficiently large $t$.
\end{prop}

\begin{proof}
Let $f$, and $f'$ and $F$ in $\RR_+[x,y,z]_n$, $\RR_+[x,y,z]_{n'}$ and $\RR_+[x,y,z]_{n+n'}$ be the polynomials $\delta(\beta(A,B,C))$, $\delta(\beta(A',B',C'))$ and $\delta(\beta(A \oplus A', B \oplus B', C \oplus C'))$. We observed earlier that $F=ff'$ and it is then clear from the definition that $N_F=N_f+N_{f'}$. 

Let $(I,J,K) \in \Delta^{n+n'}$ and let $(x,y,z) \in \RR^3$ achieve the infimum in the definition of $u_{IJK}(F)$. Consider any $(i,j,k) \in \Delta^n$ and $(i',j',k') \Delta^{n'}$ with $I=i+i'$, $J=j+j'$ and $K=k+k'$.  We have 
\begin{multline*}
u_{IJK}(F)=N_F(x,y,z)-Ix-Jy-Kz= \\ \left( N_f(x,y,z) - ix-jy-kz \right) + \left( N_{f'}(x,y,z) - i'x-j'y-k'z \right) \geq u_{ijk}(f)+u_{i'j'k'}(f').
\end{multline*}

Taking the maximum over all $(i,j,k)$ and $(i',j,k')$, we have the result.

The proof of the last statement is omitted. 
\end{proof}

\begin{prop}
Let $\partial : \HH \to (\RR_{\geq}^n)^3$ and $\partial_0 : \HIVE_n \to (\RR_{\geq}^n)^3$ be the maps deinfed in the first section. Let $F \in \HH$. Then $\partial(F)= (1/2) \partial_0 (u(F))$.
\end{prop}

\begin{proof}
By Proposition 3.2 of \cite{FPT}, $u_{ij0}(F)$ only depends on $F(1,u,0)$. Let the roots of $F(1,u,0)$, which we know to be negative, be $-r_1^2$, \dots, $-r_n^2$,  with $r_1 \geq r_2 \geq \ldots \geq r_n$.

By Jensen's theorem, 
$$N_{F(1,u,0)}(s)=\log F_{00n}+\sum_{r_i^2 < e^s} (s-\log r_i^2).$$
(The sum is over all $i$ such that $r_i <e^s$. Then $u_{j}(F(1,u,0))= \log F_{00n}- 2 \sum_{m=j}^n \log r_i$. The result is now obvious.
\end{proof}

\subsection{The Existence of Limits}

Let $\phi$ denote the composite map $\OO^{\add} \to \RR^{\Delta}$ by $(1/2) \log \circ \delta \circ \beta \circ \alpha$.

We still have not actually constructed a map from $\OO^{\add}$ or $\OO^{\mult}$ to $\HIVE$. By Theorem~\ref{Real}, if $\gamma(t)$ is any path in $\HH$ such that $\lim_{t \to \infty} (1/t) \gamma(t)$ exists, this limit is in $\HIVE$. The following would be the most elegant way a map $\OO^{\add} \to \HIVE$ could be constructed.

\begin{question}
Let $\gamma(t)$ be a path in $\OO^{\add}$ such that $\lim_{t \to \infty} \gamma(t)/t$ exists. Does the limit $\lim_{t \to \infty} \phi(\gamma(t))/t$ necessarily exist? Is it dependent only on $\lim_{t \to \infty} \gamma(t)/t$ and not on the choice of path achieving this limit?
\end{question}

I conjecture that the answer to the first question is yes, at least for ``nice'' paths. On grounds of elegance, the second statement should be true, but I find it hard to imagine how it could occur, as the map $\delta$ depends not only on the asymptopics of the entries in the matrices $X$, $Y$ and $Z$ but on those of all their minors. A weaker conjecture is that $\lim_{t \to \infty} (1/t) \phi(tA,tB,tC)$ exists for all $(A,B,C) \in \OO^{\add}$.  The primary difficulty lies in the map $\alpha$. Because of the analytic nature of its definition, it is difficult to find any data from which to extrapolate.
 
If such a limiting map does exist then, by Theorem~\ref{Torii}, it will be a map from a manifold to a polyhedral cone whose fibers are tori. If $(\a, \b, \c)$ are fixed, the fibers will be of the same dimension as the image. This suggests that in some way, we are seeing a degeneration of $\OO^{\add}$ to a toric variety.

\subsection{More Matrices}

One can ask to investigate the space $\{ (A_1, \ldots, A_r) \in \Her : \sum A_i=0 \} / U$, for $r=3$ this is $\OO^{\add}$. One can define the map $\phi$ analogously to before. One is now interested in characterizing the space $\HH^r \subset \RR_+[x_1, \ldots, x_r]_n$ of all polynomials of the form $\det(\sum x_i X_i)$ with $X_i \in \PD$. We will call a polynomial in $\HH^r$ a \emph{Vinnikov hypersurface}. One can also hope to study $\HIVE^r$ which we \emph{define to be} $l(\HH^r(\TCC))$. Write $\Delta^r=\{ (i_1, \ldots, i_r) \in \ZZ_{\geq 0}^r : \sum i_j=n \}$. So $\HIVE^r \in \RR^{\Delta^r}$. 

The first difficulty in this problem is that Theorem~\ref{AnyCurve} is not true for $r>3$. This is easy to see: giving $r$ matrices in $\Her_n$ involves $rn^2$ parameters while giving a polynomial in $\CC[x_1, \ldots, x_r]_n=\CC^{\Delta^r}$ involves $\binom{n+r-1}{r-1}$ parameters. For $r>3$ and $n$ large, $\binom{n+r-1}{r-1} > rn^2$.

The if and only if in Vinnikov's criterion, therefore, can no longer hold, as all of the topological conditions are preserved under perturbing the polynomial $F$. However, one can still derive topological properties of Vinnikov hypersurfaces.

\begin{prop}
Let $X_1$, \dots, $X_r \in PD_n$ and set $F(x_1, \ldots, x_r)=\det \left( \sum x_i X_i \right)$. Then any line meeting $\Tau_{++ \cdots +} $ meets $Z(F)$ $n$ times. If $Z(F)$ is smooth, it consists of $\lfloor n/2 \rfloor$ components that divide $\RR\PP^{r-1}$ into two pieces, with $\Tau_{++ \cdots +}$ on the contractable portion and possibly one more piece that does not disconnect $\RR\PP^{r-1}$.\end{prop}

\begin{proof}
This is analogous to the remark after the statement of Vinnikov's Criterion.
\end{proof}

We will say that any hypersurface in $\RR \PP^{r-1}$ with the above topology is a \emph{topological Vinnikov hypersurface}. It is easy to prove that the derivative of a polynomial cutting out a topological Vinnikov hypersurface itself defines a topological Vinnikov hypersurface, which allows us to prove some easy inequalities by induction as in the proof of theorem~\ref{Backward}. I do not know whether the derivative of a Vinnikov hypersurface is a Vinnikov hypersurface.

I have worked out the case of $\HIVE^4_2$. 
\begin{prop}\label{HIVE4}
$\HIVE^4_2$ is defined by the three symmetric permutations of the following inequality:
$$h_{1100}+h_{0011} \leq \max( h_{1010}+h_{0101}, h_{1001}+h_{0110} )$$
and the twelve permutations of the following inequality
$$h_{2000}+h_{0110} \leq h_{1100}+h_{1010}.$$
\end{prop}

Note that $\HIVE^4_2$ is a fan, not a cone.

\end{document}